\documentclass[
 reprint,
 amsmath,amssymb,
 aps,
]{revtex4-2}

\usepackage{hyperref}
\usepackage{xcolor}
\usepackage{graphicx}
\usepackage{dcolumn}
\usepackage{bm}
\usepackage{overpic,pict2e} 
\usepackage{tikz} 
\usetikzlibrary{positioning}

\begin{document}

\preprint{APS/123-QED}

\title{\textcolor{blue}{Exploring the electric field around a loop of static charge:\\ 
Rectangles, stadiums, ellipses, and knots}}

\author{Max Lipton}
\email{ml2437@cornell.edu}
\affiliation{Mathematics Department, Cornell University, Ithaca, NY 14853}%
\author{Alex Townsend}%
 \email{townsend@cornell.edu}
\affiliation{Mathematics Department, Cornell University, Ithaca, NY 14853}%
\author{Steven H. Strogatz}%
 \email{strogatz@cornell.edu}
\affiliation{Mathematics Department, Cornell University, Ithaca, NY 14853}%

\date{\today}

\begin{abstract}

\textcolor{blue}{We study the electric field around a continuous one-dimensional loop of static charge, under the assumption that the charge is distributed uniformly along the loop. For rectangular or stadium-shaped loops in the plane, we find that the electric field can undergo a symmetry-breaking pitchfork bifurcation as the loop is elongated; the field can have either one or three zeros, depending on the loop's aspect ratio. For knotted charge distributions in three-dimensional space, we compute the electric field numerically and compare our results to previously published theoretical bounds on the number of equilibrium points around charged knots. Our computations reveal that the previous bounds are far from sharp. The numerics also suggest conjectures for the actual minimum number of equilibrium points for all charged knots with five or fewer crossings. In addition, we provide the first images of the equipotential surfaces around charged knots, and visualize their topological transitions as the level of the potential is varied.}

\end{abstract}

\maketitle


\section{\label{sec:intro}Introduction}

\textcolor{blue}{In a first course on electricity and magnetism, students are often asked to solve problems in electrostatics~\cite{purcell1965electricity}. For example, calculating the electric field around an infinite line of uniformly distributed static charge provides good practice in working with Coulomb’s law or Gauss’s law. Another classic exercise is to calculate the electric field at all points on the symmetry axis above the center of a uniformly charged circular ring.} 

\textcolor{blue}{Here we explore the electrostatics of one-dimensional charge distributions that are less standard than lines and circles. First, working with charged loops confined to a plane, we show that charged rectangles and stadiums can give rise to symmetry-breaking bifurcations in their surrounding electric fields. Specifically, imagine stretching a rectangular charge distribution along its major axis. When the charged rectangle becomes sufficiently elongated, the electric field it generates can have three equilibria: one at the center of the rectangle, and another two located on the major axis, symmetrically placed on opposite sides of the center. These two additional equilibria emerge from a supercritical pitchfork bifurcation at a critical aspect ratio that can be calculated explicitly. Surprisingly, nothing like this happens for a charged ellipse; at all aspect ratios, its surrounding electric field has only one equilibrium point at its center.}

\textcolor{blue}{Next we allow our charged loops to move out of the plane and into three dimensions. In particular we consider charged trefoil knots, figure-eight knots, and other charged loops that have knots tied in them. Charged knots may sound at first like a contrived thing to study, but they  actually arise quite commonly in nature; they are found in long closed polymers and specifically in long loops of DNA~\cite{dommersnes2002knots, arsuaga2002knotting, weber2006numerical, orlandini2017dna}. Many new questions about equilibrium points and equipotential surfaces arise in this context that we hope will appeal to the general physics community.} 

\textcolor{blue}{Before we proceed with the analysis, let us clarify how our work differs from that recently presented elsewhere by the first author~\cite{lipton2020morse, lipton2021bound}. The articles~\cite{lipton2020morse, lipton2021bound} are aimed at specialists in topology, geometry, and knot theory, whereas the present treatment is intended for physicists. We are also more concerned here with giving physical, visual, and numerical results and examples, whereas the earlier articles concentrated on the precise statements and proofs of certain theorems we cite below. For example, it was proven in Ref.~\cite{lipton2021bound} that the electric field around any uniformly charged knot must have at least $2t + 1$ equilibrium points, where $t$ is a topological invariant known as the knot’s tunnel number. We now show numerically that this lower bound is not sharp, and we offer conjectures for the actual minimum number of equilibrium points for each type of knot with five or fewer crossings. Using computer graphics, we also visualize the equipotential surfaces of charged knots for the first time. Furthermore, our results about charged rectangles, stadiums, and ellipses are new.} 

\textcolor{blue}{\section{\label{sec:planar}Charged Loops in the Plane: Rectangles, Stadiums, and Ellipses}}

\begin{figure} 
\centering
\includegraphics[width=.35\textwidth,trim=140 560 300 120,clip]{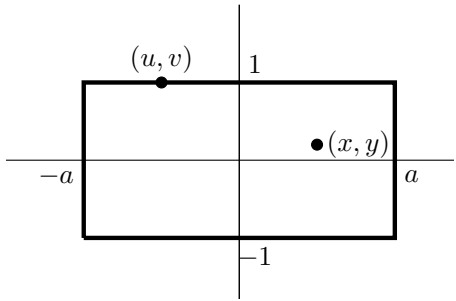}
\caption{A uniform charge distribution in the shape of a rectangle. Static charge is confined to this one-dimensional loop, and we seek to understand the electric field that surrounds it in three dimensions.}
\label{fig:ChargedRectangle} 
\end{figure}

\textcolor{blue}{\subsection{\label{sec:rectangle}Rectangle}
As a first example, consider a one-dimensional loop of static charge in the $xy$ plane, uniformly and continuously distributed in the shape of a rectangle with aspect ratio $a \ge 1$. As shown in~\autoref{fig:ChargedRectangle}, the rectangle consists of two horizontal line segments extending from $x=-a$ to $a$, lying one unit above and below the $x$ axis, capped off by vertical lines at both ends.} 

\textcolor{blue}{One of the main topics of interest to us is the location of the zeros of the electric field in the three-dimensional space around the rectangle. These zeros correspond to equilibrium points where a test charge could remain at rest (though not stably, as a well-known theorem forbids the existence of stable equilibria in an electrostatic field~\cite{purcell1965electricity}). By symmetry, there must be a zero at the center of the rectangle. But can there be any other zeros?}

\textcolor{blue}{Clearly there cannot be any zeros above or below the $xy$ plane, because all the charges in the rectangle would attract or repel a test charge in the same direction relative to the $z$-axis.  So it suffices to look for zeros in the $xy$ plane. These correspond to critical points of the potential, where the gradient vanishes. To find such zeros, if they exist, let us calculate the potential on the $xy$ plane generated by the rectangular charge distribution.} 

\textcolor{blue}{Consider the contribution to the electric potential from the top side of the rectangle. Assume a linear charge density of unit strength, for simplicity. As shown in~\autoref{fig:ChargedRectangle}, the distance between a typical line element $dx$ at $(u,v)$ on the top and a point $(x,y)$ inside the rectangle is $\sqrt{(u-x)^2+(1-y)^2}$, since $v=1$ on the top of the rectangle. For the inverse-square electrical force of Coulomb's law in three dimensions, the contribution to the potential is inversely proportional to this distance. Hence, by integrating over all the line elements on the top of the rectangle, we obtain a contribution to the potential of $$ \int_{-a}^a \frac{dx}{\sqrt{(u-x)^2+(1-y)^2}}.$$}  

\textcolor{blue}
{
\noindent This integral evaluates to 
\begin{equation}
\begin{split}
    \frac{1}{2} \log \left( \frac{a-x+\sqrt{(x-a)^2+(y-1)^2}}{a-x - \sqrt{(x-a)^2+(y-1)^2}}\right) \\
+ \frac{1}{2} \log \left( \frac{a+x+\sqrt{(x-a)^2+(y-1)^2}}{a+x - \sqrt{(x-a)^2+(y-1)^2}}\right). 
\end{split}
\end{equation}
\noindent The contributions to the potential from the other three sides of the rectangle can be evaluated similarly. In this way one can obtain the potential everywhere inside the rectangle in the $xy$ plane. (We do not show the final result because it is unpleasant to look at, but it is a sum of four terms like that above.) 
}

\begin{figure} 
\begin{minipage}{.14\textwidth}
\begin{overpic}[width=\textwidth]{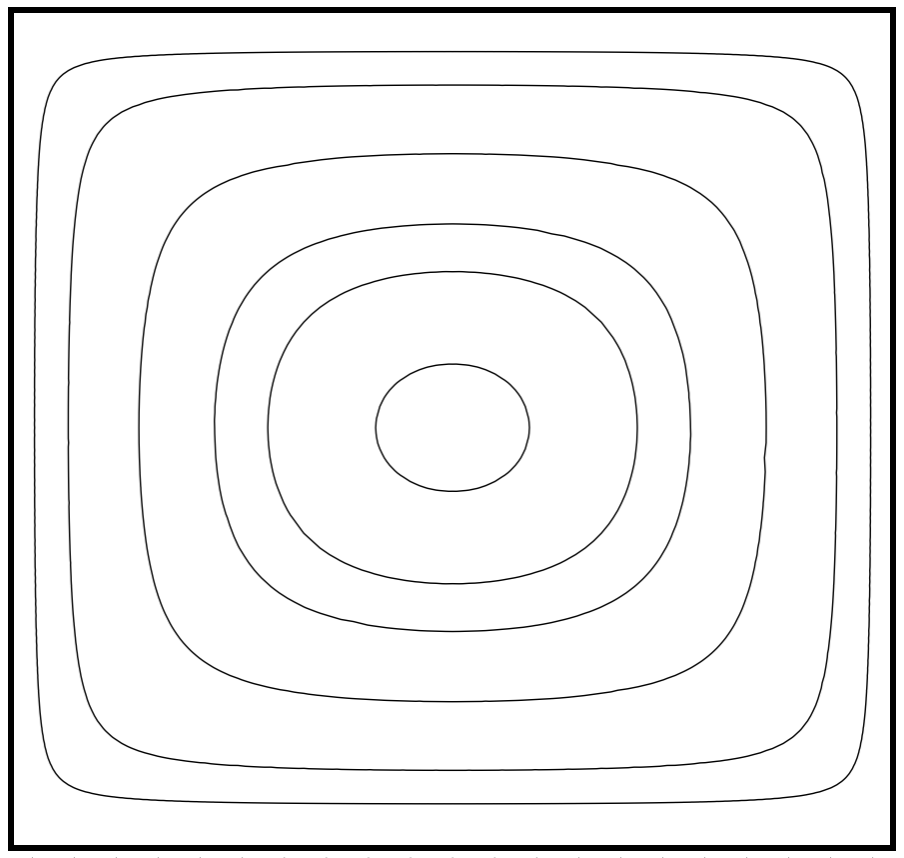}
\put(0,99) {(a)}
\end{overpic} 
\end{minipage} 
\begin{minipage}{.33\textwidth}
\begin{overpic}[width=.99\textwidth]{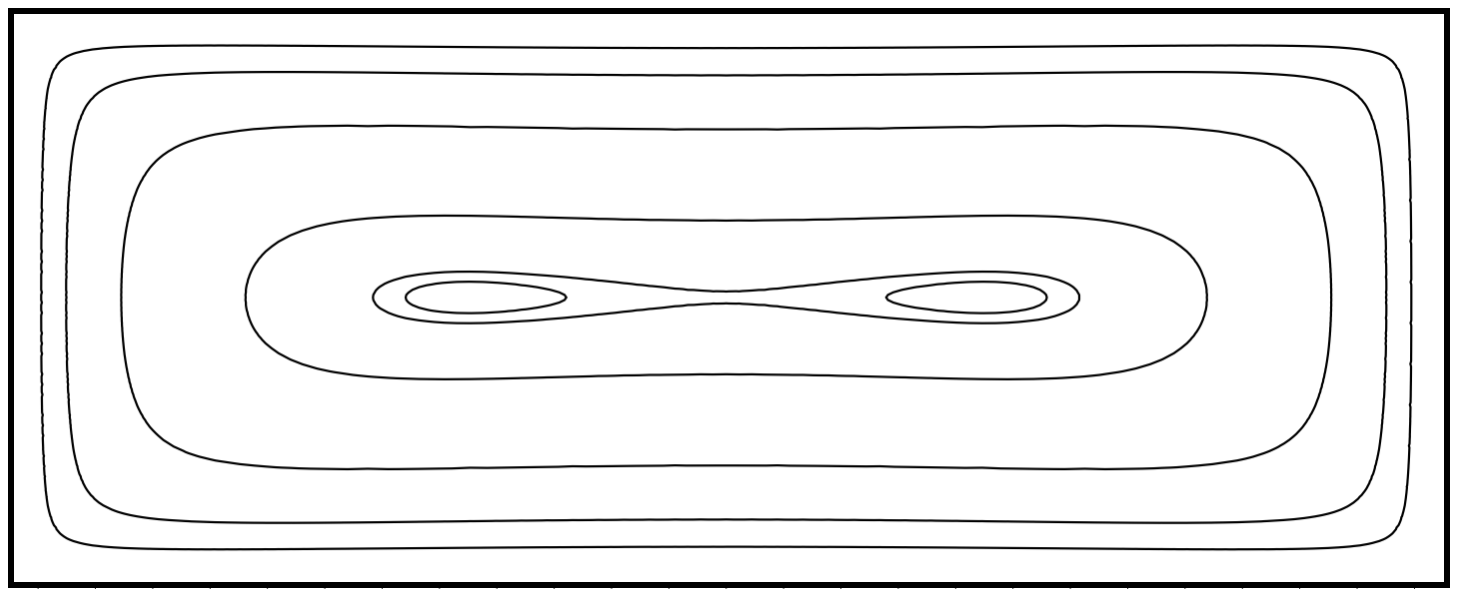} 
\put(0,42) {(b)}
\end{overpic} 
\end{minipage} 
\caption{Contour map of the electrostatic potential inside a uniform rectangular charge distribution. (a) Aspect ratio $a = 1.1$. The only critical point (also known as a zero or an equilibrium point of the electric field) occurs at the center of the rectangle. (b) Aspect ratio $a = 2.5$. There are now three zeros.}
\label{fig:RectanglePlots} 
\end{figure}

\textcolor{blue}
{
\autoref{fig:RectanglePlots} shows a contour map of the potential $\phi(x,y,z)$ restricted to the plane $z = 0$ for two values of the aspect ratio $a$. When $a$ is close to 1 and the rectangle is almost square, the only zero is at the origin (\autoref{fig:RectanglePlots}(a)). However, if we increase $a$ to 2.5, the rectangle becomes more elongated and we now find three zeros: one at the origin, and a symmetric pair on either side of the origin (\autoref{fig:RectanglePlots}(b)). 
}

\begin{figure} 
\begin{minipage}{.23\textwidth}
\begin{overpic}[width=\textwidth]{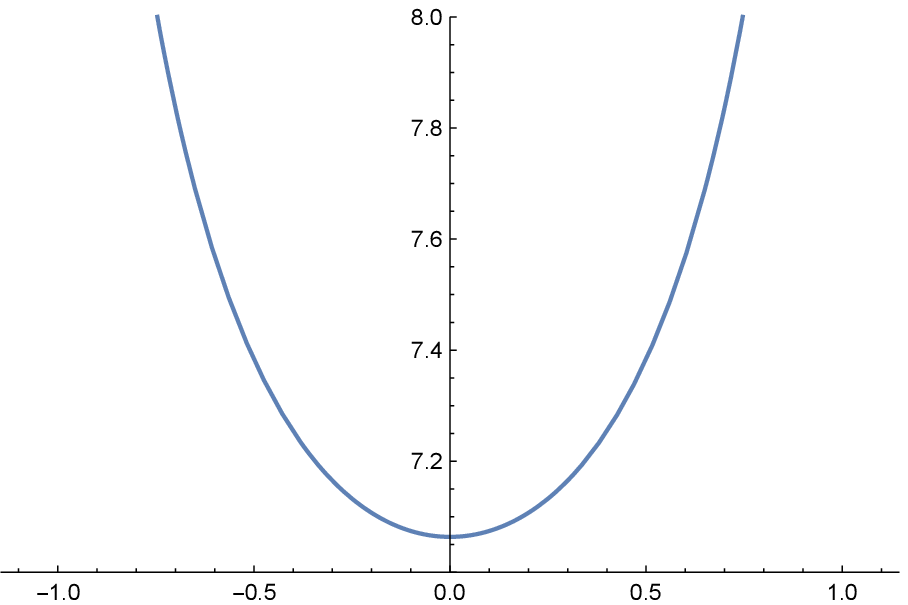}
\put(5,55) {(a)}
\put(48,-5) {$x$}
\end{overpic} 
\end{minipage} 
\begin{minipage}{.23\textwidth}
\begin{overpic}[width=\textwidth]{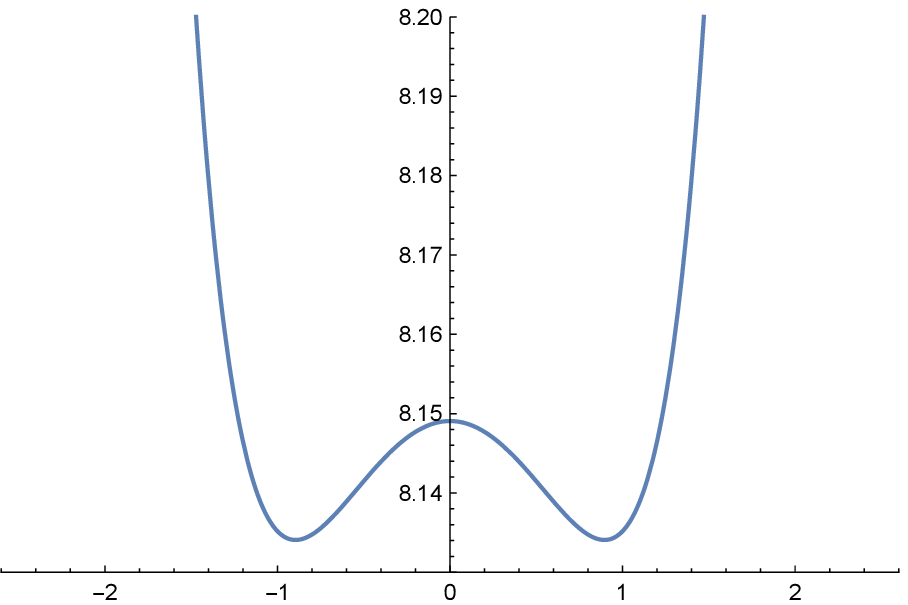}
\put(5,55) {(b)}
\put(48,-5) {$x$}
\end{overpic} 
\end{minipage}
\caption{The potential of charged rectangles along the $x$-axis. (a) Aspect ratio $a = 1.1$. The only critical point lies at the origin. (b) Aspect ratio $a = 2.5$. The three critical points where the slope vanishes correspond to electrostatic equilibrium points, i.e., zeros of the electric field.}
\label{fig:Potential1}
\end{figure}

\textcolor{blue}
{To find the threshold value of $a$ at which the bifurcation occurs, we examine how the potential varies along the $x$-axis.  \autoref{fig:Potential1} shows that the potential on the $x$-axis changes from having a minimum to a maximum at $x=0$ as $a$ increases. At the bifurcation value of $a$, the second derivative of the potential at the origin vanishes as the graph of the potential changes from concave up to concave down. By calculating this second derivative analytically, we find that the threshold value of $a$ satisfies 
$$4 + 8 a^2 -4 a^3 =0,$$
whose unique positive root is $a \approx 2.20557.$
}

\begin{figure} 
\centering
\includegraphics[width=.49\textwidth,trim=140 560 200 120,clip]{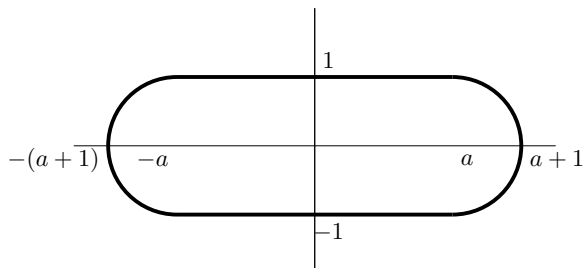}
\caption{A stadium-shaped charge distribution. The stadium's aspect ratio is defined as $a+1$, where  $a$ denotes the half-length of the straight portion of the stadium.}
\label{fig:ChargedStadium}
\end{figure}

\textcolor{blue}
{
\subsection{\label{sec:stadium}Stadium}
Our second example takes the form of a planar curve shaped like a stadium (\autoref{fig:ChargedStadium}). This curve consists of two equal parallel line segments that run from $x=-a$ to $x=a$, with $y$ values given by $y=\pm 1$. These line segments are capped off by semicircles of unit radius at either end. The resulting stadium curve is well known in studies of chaotic billiards in classical mechanics~\cite{bunimovich1979ergodic} and in experimental studies of quantum chaos~\cite{marcus1992conductance}. 
}

\textcolor{blue}
{
We now show that the electrostatic field around a uniformly charged stadium curve can have either one or three zeros, depending on the value of the stadium's aspect ratio. As before,  imagine distributing static charge uniformly and continuously along this curve in the $xy$ plane. Zeros of the electric field around it must lie in the $xy$ plane, by the reasoning given earlier for the rectangle. Indeed, the zeros are all confined to the $x$-axis, for all values of $a>0$. }

\begin{figure} 
\begin{minipage}{.19\textwidth}
\begin{overpic}[width=\textwidth]{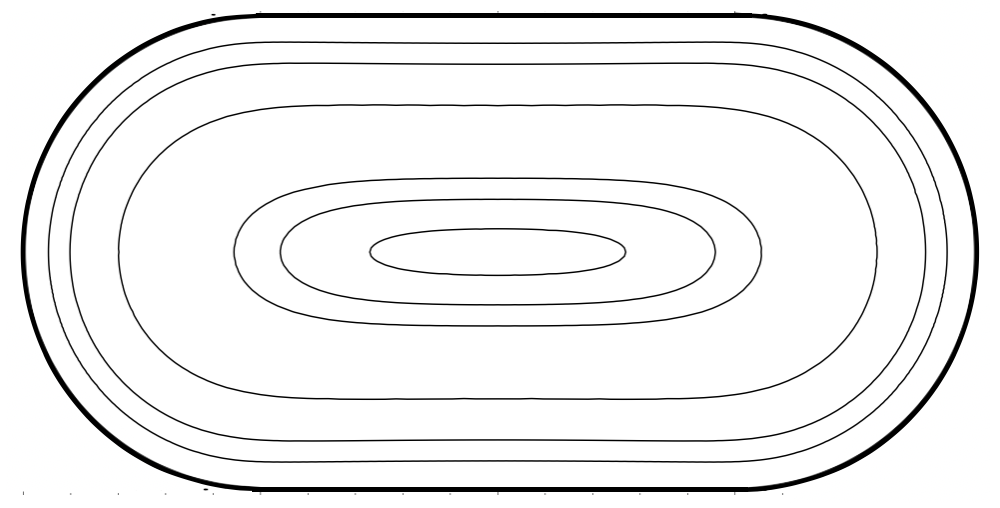} 
\put(0,46) {(a)}
\end{overpic} 
\end{minipage} 
\begin{minipage}{.28\textwidth}
\begin{overpic}[width=\textwidth]{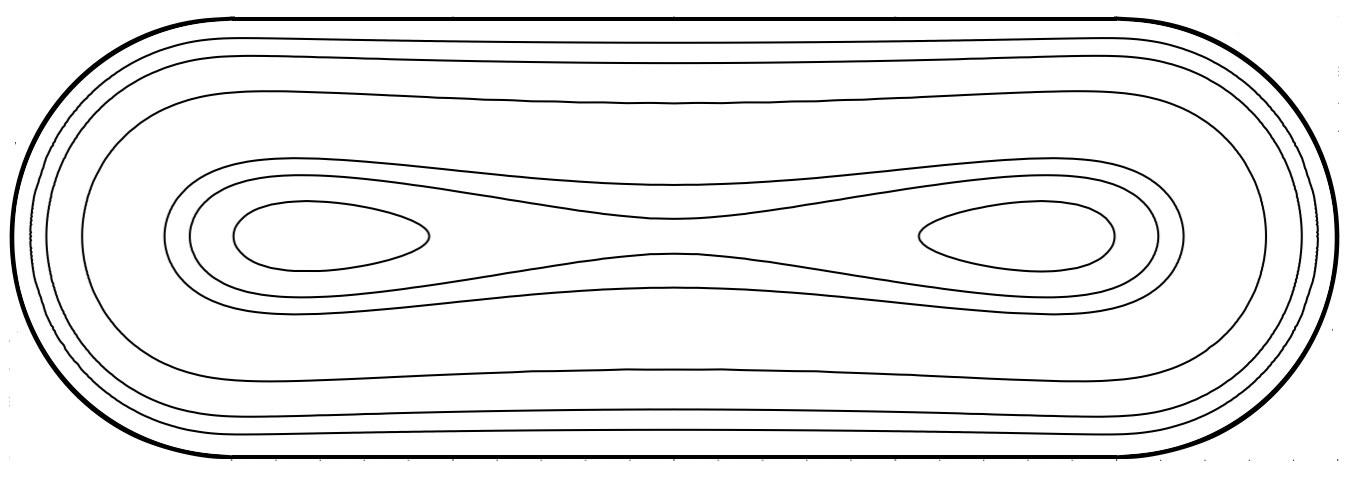} 
\put(0,33) {(b)}
\end{overpic} 
\end{minipage} 
\caption{(a) Contour map of the electrostatic potential for a charged stadium with $a = 1$. The only critical point lies at the origin. (b) Contour map for $a = 2$. Now there are three critical points: one at the origin, and a pair on either side.}
\label{fig:StadiumPotential}
\end{figure}

\textcolor{blue}
{\autoref{fig:StadiumPotential} shows the behavior of the potential for three values of $a$. For $a = 1$, there is only one zero, whereas for $a=2$ there are three zeros. 
}

\begin{figure} 
\begin{minipage}{.23\textwidth}
\begin{overpic}[width=\textwidth]{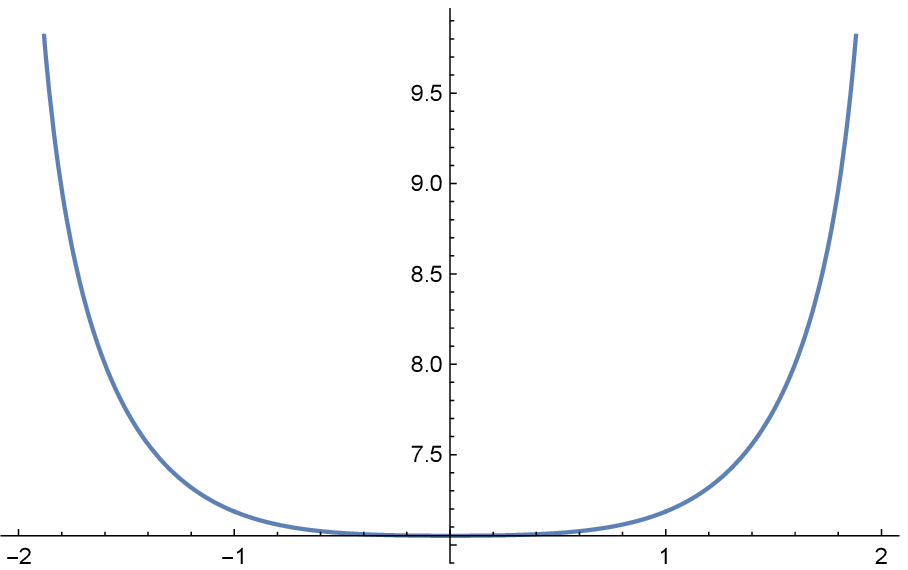} 
\put(0,62) {(a)}
\put(48,-5) {$x$}
\end{overpic} 
\end{minipage} 
\begin{minipage}{.23\textwidth}
\begin{overpic}[width=\textwidth]{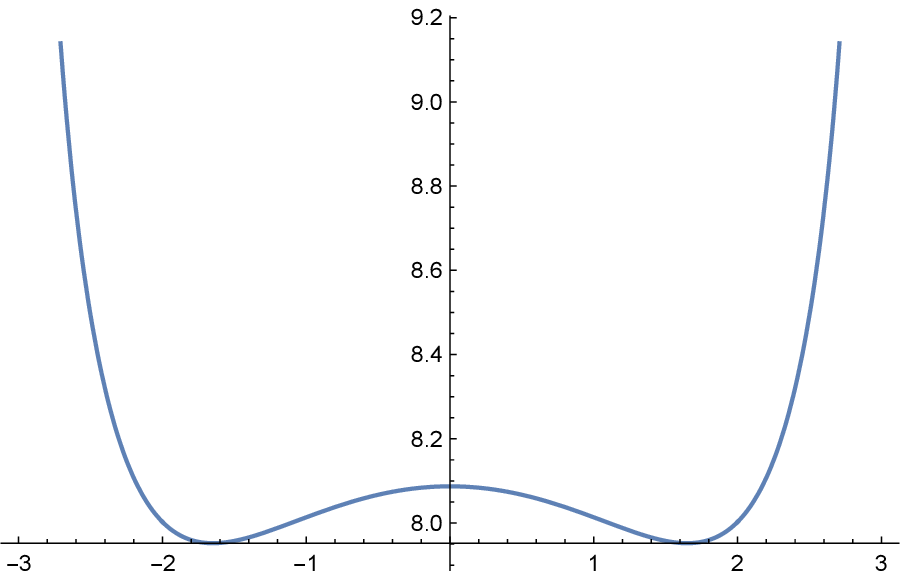} 
\put(0,62) {(b)}
\put(48,-5) {$x$}
\end{overpic} 
\end{minipage} 
\caption{(a) Potential along $x$-axis for stadium with $a = 1$. (b) Potential along $x$-axis for stadium with $a = 2$. The three critical points where the slope vanishes correspond to electrostatic equilibrium points, i.e., zeros of the electric field.}
\label{fig:StadiumPotential2} 
\end{figure} 

\textcolor{blue}
{
To find the threshold value of $a$ at which the bifurcation occurs, we follow the same procedure as  for the rectangle. The integrals arising in the expression for the potential along the $x$-axis can again be found analytically, yielding a formula in terms of elliptic functions and inverse $\sinh$ functions   (we omit the details). By solving for the value of $a$ at which the second derivative of the potential along the $x$-axis vanishes at the origin, we find the bifurcation occurs at $a \approx 1.1313$. The corresponding aspect ratio of $a+1 \approx 2.13$ is not too different from the value of $2.2$ found earlier for the rectangle.  
}

\textcolor{blue}
{
\subsection{\label{sec:ellipse}Ellipse}
}

\begin{figure} 
\centering
\begin{overpic}[width=.23\textwidth]{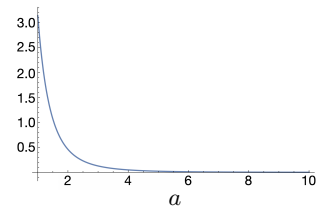} 
\put(50,-5) {$a$}
\end{overpic}
\caption{The second derivative in $x$ of $\phi$ at the origin for a charged ellipse with aspect ratio $a$. Here $a$ is defined as the length of the ellipse's major axis, measured in units of its minor axis. Numerically, it appears that the second derivative of the potential is strictly positive at the origin. This observation provides numerical evidence that the electric field associated with any charged ellipse has exactly one zero.}
\label{fig:DerivativeEllipse}
\end{figure}

\textcolor{blue}
{Curiously, when we redo the calculations above for an ellipse, we always find only a single zero (the trivial equilibrium point at the origin expected by symmetry), no matter what aspect ratio we choose. We have not proven this, but we find numerically that the second derivative of the potential (restricted to the $x$-axis) is strictly positive for all aspect ratios and hence all possible elliptical shapes (\autoref{fig:DerivativeEllipse}).
}

\textcolor{blue}
{
An intuitive explanation of this result relies on the shape difference between elongated rectangles and stadiums on the one hand, and elongated ellipses on the other. The key is that there is enough charge on the short side of a rectangle (or the small semicircle of a stadium) to counterbalance the forces exerted by the rest of the shape, as long as a test charge is situated sufficiently close to the short side (or the small semicircle). For example, consider a highly elongated rectangle. It seems plausible that there should be an equilibrium point very close to either of the short sides of the rectangle; all the charge on that side can exert enough electrical force to counterbalance the weaker forces exerted by the much more distant charges located on the rest of the rectangle. The same effect holds for a stadium. But with an ellipse, there is apparently not enough charge on the narrow end of the ellipse to offset the force coming from the rest of the shape. 
}

\section{\label{sec:knots} Charged Knots}

\textcolor{blue}{In the rest of this article we move out of the plane and consider a continuous charge distribution in the shape of a knotted loop in three-dimensional space.  Generalizing our work in the previous section, we ask:} How many equilibrium points must exist, and where are they located in relation to the knot? What do the equipotential surfaces look like, and how does their topology change as we vary the level of the potential?

\subsection{\label{sec:knot_intro} Background}

The mathematical theory of knots and links grew out of two problems in classical physics~\cite{przytycki1998classical}. In 1833, while studying electromagnetism and mutual inductance, Gauss discovered a formula for the linking number of two closed curves in three-dimensional space~\cite{ricca2011gauss}. A few decades later, Lord Kelvin proposed his vortex model of the atom, in which different elements were imagined to arise from different patterns of knotted and linked vortex rings in the ether~\cite{thomson1867vortex}. Since then, knots and links have been studied in many other parts of physics~\cite{kauffman1991knots, stasiak1998ideal, calvo2005physical, adams2021encyclopedia}, including quantum field theory~\cite{witten1988topological, witten1989quantum}, liquid crystals~\cite{irvine2014liquid}, plasmas~\cite{berger_field_1984, smiet2015self}, Bose--Einstein condensates~\cite{kawaguchi2008knots, hall2016tying}, fluids~\cite{moffatt1969degree, kleckner2013creation}, superfluids~\cite{kleckner2016superfluid}, the biophysics of DNA and other polymers~\cite{dommersnes2002knots,weber2006numerical}, the mechanics of ropes and elastic rods ~\cite{buck1999thickness, audoly2007elastic, jawed2015untangling, moulton2018stable, patil2020topological}, classical field theory~\cite{faddeev1997stable}, knotted scroll waves in excitable media~\cite{winfree1984organizing, sutcliffe2003stability}, $n$-body choreographies in celestial mechanics~\cite{calleja2019torus}, and electromagnetic waves and fields~\cite{ranada1990knotted, lomonaco1996modern, berry2001knotted, leach2004knotted, irvine2008linked, kedia2013tying, arrayas2017knots, larocque2018reconstructing}.

\subsection{\label{sec:physical_argument} Physical Intuition about Equipotential Surfaces and Equilibrium Points}

\textcolor{blue}{To get an intuitive feel for knotted charge distributions, let us start with the ``trivial knot'' (also known as the ``unknot''), and for simplicity, consider its most symmetrical realization: a perfectly circular and uniform loop of static charge.} By symmetry, or by using Coulomb's law, one can prove there is an equilibrium point at the center of the circle. Moreover, the symmetry of the situation suggests that the equipotential surfaces enclosing the circle should be nested tori, at least sufficiently close to the circle. Far from the circle, they must approach spheres, because a circular loop of charge looks like a point source in the far field.


\begin{figure*}
\includegraphics[scale=1]{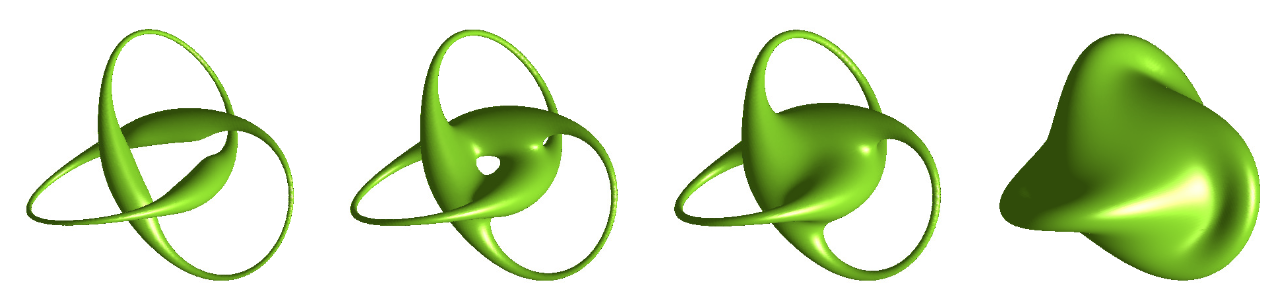}
\caption{Topologically different equipotential surfaces for a charged trefoil. Only ``regular'' values of the potential are shown, at which the equipotential surfaces are smooth manifolds. Bifurcations (not shown) lie in between these examples. In the leftmost panel, the potential is very large and positive close to the knot, and the corresponding equipotential surface resembles a thin tube around the knot. As we lower the potential (corresponding to moving from left to right in the figure), the surface inflates and swells up. At bifurcation values, it self-intersects and creates new zeros in the electric field, corresponding to new equilibrium points. The equipotential surfaces shown here have genus values of $1, 4, 3,$ and $0$, moving from left to right.}
\label{fig:trefcode}
\end{figure*}

\begin{figure*}
\begin{minipage}{.245\textwidth}
\includegraphics[width=1.4\textwidth,trim=3cm 2cm 0cm 1cm,clip]{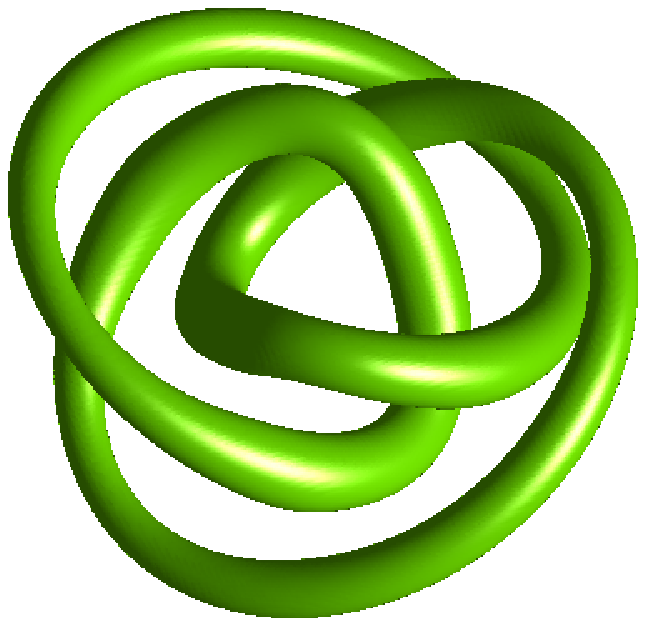}
\end{minipage} 
\begin{minipage}{.245\textwidth}
\includegraphics[width=1.4\textwidth,trim=3cm 2cm 0cm 1cm,clip]{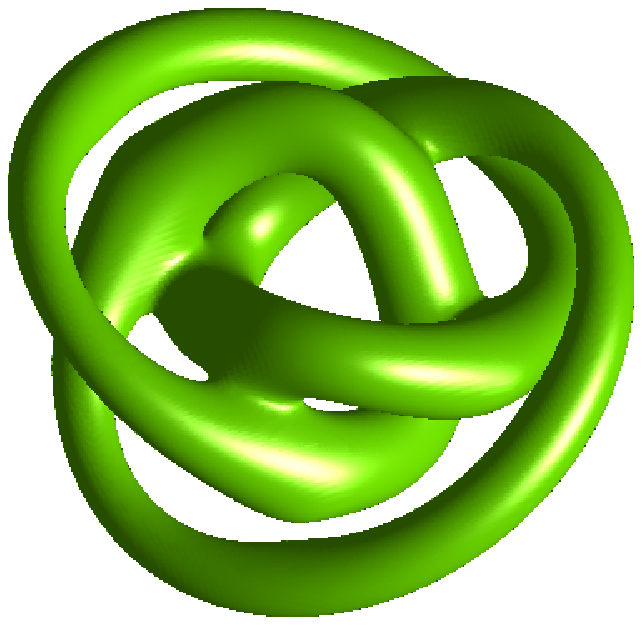}
\end{minipage} 
\begin{minipage}{.245\textwidth}
\includegraphics[width=1.4\textwidth,trim=3cm 2cm 0cm 1cm,clip]{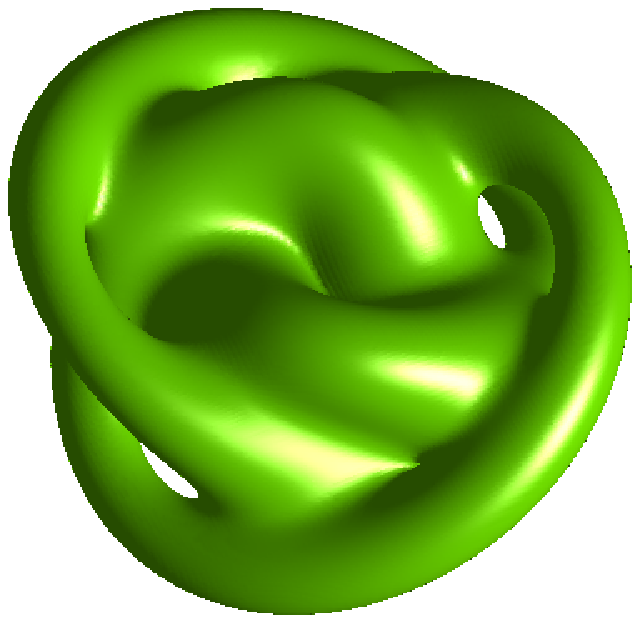}
\end{minipage} 
\begin{minipage}{.245\textwidth}
\includegraphics[width=1.4\textwidth,trim=3cm 2cm 0cm 1cm,clip]{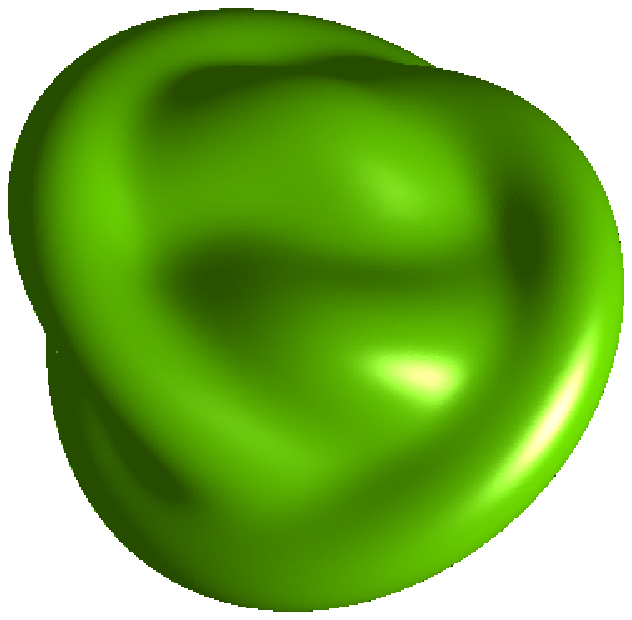}
\end{minipage} 
\caption{Topologically different equipotential surfaces for a charged figure-of-eight knot.}
\label{fig:FigureOfEight}
\end{figure*}


Remarkably, a similar line of reasoning works even if the charged loop has knots tied in it. Although it becomes impossible to calculate the location of the equilibrium points analytically in these cases, we can still prove that equilibria must exist in the surrounding electric field. 

To do so, picture the equipotential surfaces. \autoref{fig:trefcode} and~\autoref{fig:FigureOfEight} show what they look like for a trefoil knot and a figure-eight knot, but the same idea works for any sufficiently smooth knot.  One expects that close to the knot, the equipotential surfaces must be tubular versions of the knot itself---in other words, they must be knotted tori---whereas far from the knot, they must resemble spheres (because a charged knot, like any other compact charge distribution, appears point-like in the far field). So if we imagine continuously varying the potential from high levels near the knot to low levels at infinity, the equipotential surfaces must continuously deform from knotted tori into spheres. To make this transition, the knotted tori swell up, collide with distant parts of themselves, and reconnect in ways that alter their topology. At such collisions, two patches of a single equipotential surface intersect tangentially. Because the electric field vector lies along the normal to each of the colliding patches, and because those normal vectors point in opposite directions at the point of collision, the electric field must vanish there. 

This argument suggests that there must be one or more zeros in the electric field around any charged knot. Each zero represents an equilibrium point where a test charge could remain motionless. As is well known, there are no \emph{stable} equilibria in an electrostatic field, so all these equilibrium points must be unstable \cite{purcell1965electricity}. Indeed, they are all saddle points, with either one- or two-dimensional unstable manifolds. 

These intuitive ideas have recently been sharpened and made rigorous with the help of Morse theory, algebraic topology, and geometric topology~\cite{lipton2020morse, lipton2021bound}. In particular, one theorem~\cite{lipton2021bound} provides a lower bound on the number of equilibrium points around any charged knot. It states that the electric field must have at least $2t + 1$ zeros, where $t$ is a topological invariant known as the tunnel number of the knot \cite{clark1980heegaard, morimoto1996identifying, morimoto2016tunnel}. 

But as we will show below, this $2t+1$ lower bound turns out to be rather loose. For example, the tunnel number for a trefoil knot is known to be 1 (indeed, the tunnel number is 1 for all prime knots with seven or fewer crossings~\cite{morimoto1996identifying, morimoto2016tunnel}), so the $2t+1$ bound implies that a charged trefoil must always have at least three zeros in its surrounding electric field. Yet we have never seen fewer than seven zeros in our computations, no matter how the trefoil is bent, twisted, or otherwise deformed. So is seven the absolute minimum? Or could there be some needle-in-a-haystack conformation of a charged trefoil that has three, four, five, or six zeros? The problem is currently unsolved, and is just one of the many open problems about charged knots. 

In what follows, we propose conjectures for the minimum number of zeros for all prime knots with up to five crossings, based on our numerical experiments. But proving (or disproving) these conjectures---and extending them to a much wider range of knots---remains a challenge. 

One may reasonably ask why we are so concerned with these zeros, given that they represent unstable equilibria and are therefore of questionable physical significance. Although it remains to be seen whether they are physically important, from a mathematical standpoint they are definitely interesting for the following reason. The minimum number of zeros for the electric field produced by a charged knot of a given knot type provides a new knot invariant, and one based in elementary physics. It may be a difficult invariant to calculate at this early stage. But we are intrigued by the possibility that subsequent research will find connections between it and other, better established knot invariants.  

Our simulations also allow us to investigate how the equipotential surfaces change their topology as we vary the level of the potential. This sequence of topological transitions can be characterized with a set of integers we call the Morse code for the knot~\cite{lipton2020morse}. We compute these transitions numerically for some convenient embeddings of the simplest knots and illustrate them with computer graphics (\autoref{fig:trefcode} and \autoref{fig:FigureOfEight}).


\subsection{\label{sec:numerical} Numerical Methods}

To describe our computations, we introduce some notation and terminology.  Let the knot $K$ be parametrized by a vector-valued function $\pmb{r}(t)$, where $0 \le t \le 2 \pi$. Because the knot forms a closed loop, we also require that $\pmb{r}(0) = \pmb{r}(2\pi)$. 
Then, from Coulomb's law, the electric potential $\phi$ at a point $\pmb{x}\in\mathbb{R}^3$ away from the knot is given in dimensionless form by
\begin{equation}
\phi(\pmb{x}) = \int_{\pmb{r} \in K} \frac{|\pmb{dr}|}{|\pmb{x}-\pmb{r}|} = \int_0^{2\pi} \frac{|\pmb{r}'(t)| dt}{|\pmb{x}-\pmb{r}(t)|},
\label{eq:KnotPotential}
\end{equation} 
where $|\cdot|$ denotes the magnitude of a vector quantity. (We have written the potential in dimensionless form for convenience; one could include physical parameters like the vacuum permittivity $\epsilon_0$ or the uniform  charge density $\rho$ along the knot, but we have chosen not to do so as they play no role in our analysis.) The electric field associated with the potential is given by $\pmb{E} (\pmb{x}) = -\nabla \phi(\pmb{x})$. The zeros (i.e., the equilibrium points) of the electric field are equivalent to the critical points of the electric potential $\phi$; as such, we will continue to use the terms zeros, equilibrium points, and critical points interchangeably, as we did earlier in Section II. 

For most knots, the potential $\phi(\pmb{x})$ and its critical points cannot be calculated analytically. We must rely on numerical integration and rootfinding techniques. To perform these computations, we replace the continuous knot by $N+1$ unit point charges located at $\pmb{r}(t_0),\ldots, \pmb{r}(t_N)$, where $t_k = 2\pi k/(N+1)$, and use the following trapezoidal approximation of~\eqref{eq:KnotPotential}:
\begin{equation}
\phi(\pmb{x}) \approx \frac{2\pi}{N+1}\sum_{j=0}^N \frac{|\pmb{r}'(t_j)|}{|\pmb{x}-\pmb{r}(t_j)|}.
\label{eq:trapapprox}
\end{equation}
We use the same trapezoidal approximation for the electric field $\pmb{E} (\pmb{x}).$ (There are more efficient approaches for evaluating $\phi(\pmb{x})$ and $\pmb{E} (\pmb{x})$ when $N$ is very large; these are based on multipole expansions such as those used in the fast multipole method~\cite{beatson1997short}.) 

As mentioned above, we are interested in finding the zeros of the electric field. These are defined as points  $\pmb{x^*}\in\mathbb{R}^3$ such that $\pmb{E}(\pmb{x^*})= \pmb{0}$. To find them, we start with initial guesses and then refine the guesses iteratively, using a multivariable Newton method.  To obtain reasonable initial guesses, we use an algorithm known as ``3D Marching Cubes,'' a computer graphics algorithm for finding level sets of a scalar function~\cite{lorensen1987marching}. Here, we use Marching Cubes on each of the three components of $\pmb{E}(\pmb{x})$ to find their zero level sets. The algorithm partitions a large cuboid containing the knot into small cubes; then, on each cube, it uses a bilinear approximation of that component. If the bilinear approximations of all three components of $\pmb{E}$ pass through zero in the same cube, then we take the center of that cube as an initial guess for our multivariable Newton method. Some initial guesses to Newton end up diverging or converging to a far away critical point, and we throw these away. In contrast, the successful initial guesses quickly converge to the approximate locations of the zeros of $\pmb{E}(\pmb{x})$.

Along with the zeros, we are also interested in the equipotential surfaces. These are given by $\phi^{-1}(v)$, where $0 < v < \infty$ is some given voltage level. The relevant values of $v$ range from small positive values far from the knot, to large positive values close to the knot. Let $v^{*} = \phi(x^{*})$ denote the potential at an equilibrium point. Recall that equipotential surfaces undergo self-collisions at $x^{*}$ and lose smoothness there. So to get a smooth surface, we perturb $v^{*}$ to a nearby regular value $v$ at which the Hessian matrix of second derivatives of $\phi$ (equivalent to the Jacobian of $\pmb{E}$) has full rank. By the implicit function theorem, this full rank condition ensures that  $\phi^{-1}(v)$ is a smooth, orientable, compact surface without boundary. We then use the Marching Cubes algorithm to render the surface. By repeating this process for a range of $v$ values, we can explore how the equipotential surfaces change as we vary the level of the potential.

\subsection{\label{sec:knot_results} Results and Conjectures for Charged Knots}

To illustrate the results obtained with this approach, consider the following parametrization of a trefoil knot:
\noindent
\begin{equation}
\pmb{r}(t) = (\sin{t} + 2 \sin{2t}, \cos{t} - 2\cos{2t}, -\sin{3t}).
\label{trefoil_parametrization}
\end{equation}
In our numerical simulations, we sampled a cubic domain of $30 \times 30 \times 30$ initial guesses in a mesh surrounding the knot and ran the multivariable Newton method to test for convergence to a zero.  We rejected iterations that grew too large, or were within a small distance threshold from another computed zero, indicating a duplicate. 

\begin{figure}
    \centering
    \includegraphics[scale=1]{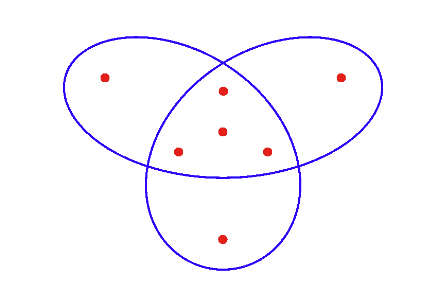}
    \caption{A top-down view of the locations of the zeros (red dots) in the electric field around a uniformly charged trefoil knot. 
    }
    \label{fig:trefcritset}
\end{figure}

\autoref{fig:trefcritset} shows that the electric field has seven zeros for this particular embedding of a trefoil. By calculating the eigenvalues at these zeros, we can confirm that they are all saddle points and classify them by their indices (the dimensions of their stable manifolds). 

Then, to obtain representatives of the equipotential surfaces around the trefoil, we compute the critical values $v^{*}$ at the zeros, perturb them to nearby values $v$, and take their inverse images $\phi^{-1}(v)$. For the parametrization~\eqref{trefoil_parametrization}, we find  the outer triplet of zeros in \autoref{fig:trefcritset} has $v^{*} \approx 12.79$ and indices of 1; the inner triplet has $v^{*} \approx 15.82$ and indices of $2$; and the origin has $v^{*} \approx 15.42$ and index $1$. To see what the equipotential surfaces look like in between these critical cases, we compute the level sets $\phi^{-1}(v)$ for the perturbed values $v = 12.7, 15, 15.5,$ and $16$. 
\autoref{fig:trefcode} shows the resulting surfaces. Topologically, they are knotted tori with various numbers of holes. 

\begin{figure} 
\centering
\includegraphics[scale=1.1]{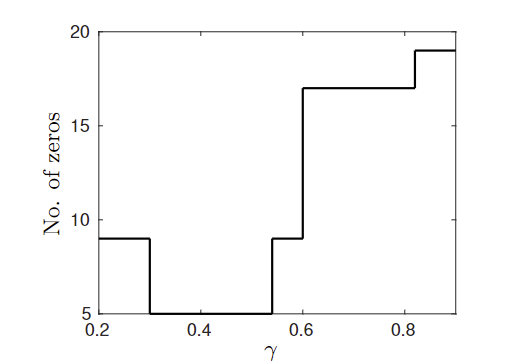}
\caption{The number of zeros in the electric field around a  figure-eight knot as it is flattened. Here, the knot is parametrized as $x(t) = (2 + \cos 2t)\cos 3t$, $y(t) =(2 + \cos2t) \sin 3t$, $z(t)=\gamma\sin4t$ and we vary the height $\gamma$.} 
\label{fig:staircase}
\end{figure}

In the example above, we assumed a highly symmetrical parametrization of a trefoil. What happens if we break the symmetry or, more generally, if we deform a knot continuously without allowing it to pass through itself? How does that affect the number of zeros in the electric field around the knot? Ideally, we would like to find deformations that cause as many zeros as possible to coalesce, thus bringing us closer to the absolute minimum number of zeros, whatever that may be, and perhaps allowing us to improve on the tunnel number bound. 

One strategy is to deform the knot so as to reduce its complexity in some way, as quantified by an energy functional or a more general Lyapunov function. Two energy functionals in the literature on physical knot theory are the M\"obius energy \cite{freedman94} and the Buck--Orloff energy \cite{buck95}, whose locally minimal configurations enjoy nice regularization properties. But we have found  a simpler strategy to be useful: we slowly squash the knot from the top down and watch what happens to its zeros. 

\autoref{fig:staircase} shows that by flattening a figure-eight knot we can reduce its number of zeros from nineteen to five; then, unexpectedly, the number goes back up to nine. It is intriguing that the number of zeros can either increase or decrease as one flattens the knot.  The staircase structure of the graph reveals that zeros can appear or disappear in pairs (which is what one expects generically~\cite{lipton2020morse}) or in two simultaneous pairs (due to non-generic symmetries in the particular parametrization of the knot). We suspect that five is the smallest number of zeros possible for \emph{any} parametrization of a  figure-eight knot, symmetrical or otherwise.

Let $Z$ denote the minimum number of zeros around a charged knot $K$, where the minimum is taken over all smooth embeddings of $K$. Although we have not found a formula for $Z(K)$, we conjecture that it can be bounded from above and below by two standard topological invariants of $K$: 
\noindent
\begin{equation}
2t+1 \le Z \le 2c+1.
\label{Z_and_bounds}
\end{equation}
We have already met the $2t+1$ lower bound, which uses the  knot's tunnel number $t$. This lower bound has been proven~\cite{lipton2021bound}. By contrast, the $2c+1$ upper bound is conjectural. It involves the knot's crossing number $c$, defined as the minimal number of crossings possible in a planar projection of the knot. We have numerical evidence and a plausibility argument but not a proof, for reasons we will explain momentarily. Table I lists these bounds along with our conjectured values of $Z$ for some simple knots. 

\begin{table}[ht]
\begin{tabular}{c c c c}
Knot & $2t+1$ & Conjectured $Z$ & $2c+1$ \\
\hline 
Unknot & 1 & 1 & 1 \\
Trefoil &  3 & 7 & 7 \\
Figure-eight & 3 & 5 & 9\\
Cinquefoil & 3 & 11 & 11\cr
Three-twist & 3 & 11 & 11 \cr 
\hline
\end{tabular}
\caption{Conjectured minimum number of zeros $Z$ for knots with five or fewer crossings. In each case, our conjectured value of $Z$ lies between the proven lower bound $2t+1$ and the conjectured upper bound $2c+1$. Our conjectures are based on numerical experiments that used the following knot parametrizations: (1) Unknot: $x(t) = \cos t$, $y(t) = \sin t$, $z(t)=0$. Trefoil: $x(t) = \sin t + 2\sin 2t$, $y(t) = \cos t - 2\sin 2t$, $z(t) = -\gamma\sin 3t$. Figure-eight: $x(t) = (2 + \cos 2t)(\cos 3t)$, $y(t) =(2 + \cos 2t) (\sin 3t)$, $z(t)=\gamma\sin 4t$. Cinquefoil: $x(t) = (\cos 2t) (3 + \cos5t)/2$, $y(t) = (\sin2t) (2 + \cos5t)/2$, $z(t)=(\gamma\sin 5t)/2$. Three-twist: $x(t) = 2\cos(2t+ 0.2)$, $y(t) = 2\cos(3t + 0.7)$, $z(t) = \gamma\cos7t$.
The standard parametrization has height $\gamma = 1$. To find a plausible conjecture for the minimum number of zeros, we flatten the knot by slowly decreasing $\gamma$ to $0$, and watch the zeros coalesce.} 
\label{Z_table}
\end{table}

Our plausibility argument for $2c+1$ as an upper bound on $Z$ is based on the following observation: For any knot $K$, we can always construct an embedding of $K$ that has \emph{at least} $2c+1$ zeros in its electric field. Unfortunately, that statement is weaker than what we need. To prove that $2c+1$ is an upper bound on $Z$, we would need to construct an embedding with \emph{exactly} $2c+1$ zeros. We suspect that this stronger construction is always possible too, but we have not proven that yet. For now, we outline the main ideas behind our weaker construction~\cite{lipton2020morse}, which proceeds in three steps. 

First, we flatten $K$ into a plane curve. The resulting curve crosses itself and bounds a number of planar regions that we call holes (more properly, ``holes'' are bounded connected components in the planar complement of the flattened knot). The electric field in $\mathbb{R}^3$ produced by a charged planar curve has a special property: 
at all points in the plane away from the curve, the electric field vector $\pmb{E}$ lies within the same plane. Moreover, at points very close to the charged curve, the direction of $\pmb{E}$ is nearly perpendicular to the curve (because it receives its dominant contribution from the portion of the charged curve nearby). Thus the winding number of the electric field around the boundary of each hole is $1$. The Poincar\'e--Hopf index theorem then implies that each hole contains \emph{at least} one source or sink zero of the planar vector field. These planar zeros extend to zeros of the full electric field in $\mathbb{R}^3$, because the out-of-plane component of $\pmb{E}$ is also zero, as discussed above.  
In Ref.~\cite{lipton2020morse} it is proven that these zeros are all saddle points of index $1$ in $\mathbb{R}^3$. A counting argument then shows that a planar curve with $c$ crossings has $c+1$ holes, and since we just showed that each hole must contain at least one zero, we arrive at our first conclusion: the electric field around a flattened knot has at least $c+1$ zeros. 

But we are not done yet. A flattened knot is not an admissible  knot because it has self-crossings. So the second step is to perturb the flattened knot by lifting one of its strands up out of the plane, ever so slightly, at each crossing to restore the topology of the original knot. By performing these lifting operations in tiny neighborhoods of the crossings that are sufficiently far away from the aforementioned zeros, we are sure to preserve the existence and topological types of the $c+1$ (or more) zeros deduced in the first step, thanks to  the structural stability of gradient vector fields \cite{palis83}. 

Now comes the third and final step. By applying the Morse inequalities \cite{nico}, one can show~\cite{lipton2020morse} that each lifting performed in the second step gives rise to a new zero at the associated crossing, and this zero is of index 2. Since a total of $c$ liftings are performed, altogether we get $c$ additional zeros of index $2$. Counting the $c+1$ (or more) zeros of index $1$ and the $c$ zeros of index $2$, we conclude the electric field around a squashed but not strictly planar version of $K$ has at least $2c + 1$ zeros. 

The weakness of the ``at least'' part of the conclusion can be traced back to the Poincar\'e--Hopf index theorem; that is where the first ``at least'' qualifier popped up. If we could ensure \emph{exactly} one zero in each hole, we would then be able to claim what we want: $Z \le 2c+1$. We suspect the uniqueness of the zero in each hole would follow if the holes were round enough (neither too elongated nor too non-convex), but this is what remains to be properly formulated and proven in future work.

\section{Discussion}

As we have tried to show in this article, the study of knotted charge distributions opens up many new directions for exploration. The questions are mainly motivated by their conceptual simplicity and theoretical appeal, but they could have real-world implications. For instance, given that the zeros of electric and magnetic fields are relevant to problems of plasma confinement in nuclear fusion and to trapping of cold atoms, related questions may be of experimental interest in those settings~\cite{hirsch1967inertial,rider1995general,hedditch2015fusion}. Charged knots also arise in molecular biology and polymer physics, where researchers study the knottedness of charged DNA molecules and their interactions with  electric fields~\cite{dommersnes2002knots, arsuaga2002knotting,weber2006numerical,orlandini2017dna}. 


\begin{acknowledgments}
We thank Greg Buck for helpful discussions. M.L. was supported in part by an NSF RTG Grant DMS-1645643. A.T. was supported by National Science Foundation grants DMS-1818757, DMS-1952757, and DMS-2045646.
\end{acknowledgments}

\bibliography{references}


\end{document}